\documentclass{amsart}

\usepackage{amscd}
\usepackage{amssymb}
\newtheorem{theorem}{Theorem}[section]
\newtheorem{lemma}[theorem]{Lemma}
\newtheorem{propo}[theorem]{Proposition}
\newtheorem{corol}[theorem]{Corollary}

\DeclareMathOperator{\sgn}{sgn}
\begin{document}

\title{Canonical Infinitesimal Deformations }
\author{Ziv Ran}
\address{}
\email{}

\maketitle

\tableofcontents

\setcounter{section}{-1}
\section{Introduction}

  Our purpose here is to develop a 'canonical' approach to infinitesimal and
formal deformation theory. For simplicity we shall stick in the paper mainly
to one fundamental-and somewhat typical-case, that of a compact complex
manifold $X$ without global vector fields. Our starting point, and model
, is the classical (first-order) Kodaira-Spencer formalism: this associates
to any deformation ${\mathfrak X}/S$ with special fiber $X$ a 'Kodaira-Spencer' map
$$\kappa : T_0S \to H^1 (X,T)$$
where $T=T_X$ is the tangent sheaf, and consequently obtains a canonical
identification between the set of first-order deformations of $X$ and the
cohomology group $H^1 (X,T)$. It is then natural to seek a higher-order analogue
of this, for n-th order tangent spaces and n-th order deformations. At a
minimum, one would like an n-th order analogue of $\kappa$:
$$\kappa_n = \kappa_n({\mathfrak X}/S) : T^{(n)}S \to (?)_n$$
where $(?)_n$ is an explicit ( and preferably computable) cohomological functor
of $X$ which, at least in favorable cases (e.g. when a global moduli space
exists), should be canonically identifiable with the n-th order tangent
space at a smooth point of  the moduli space. Put another way, one knows, when $H^0(T)=0$,
the {\it existence} of an universal formal deformation
$$\hat{X}/\hat{S}= \lim_n X_n/S_n .$$
The problem is to write each $X_n/S_n$, i.e, the universal n-th order deformation,
as an explicit cohomological functor of $X$, extending the above Kodaira-Spencer
identification of first-order deformations.

Our approach to this is to combine some earlier constructions from $[R_1]$
with an important and very novel insight coming out of recent work of
Beilinson, Drinfeld and Ginzburg, cf. [BG].  The latter is, among other
things, concerned specifically with deformations of vector bundles and
principal bundles on a fixed complex curve $X$.  It gives a formula for the
$n$-th order cotangent space to the moduli of such bundles as $H^0$ of a
suitable sheaf on the Knudsen-Mumford space $\hat{X}^n$ parametrizing $n$-tuples of points on $X$.  This suggests the simple but stunning--to the author--and
very broad philosophy that $n$-th order deformations should be related to
$n$-tuples:  e.g. that $T^{(n)} M$ or something similar ought to be
writable in terms of cohomology, at least on some sort of parameter
space for $n$-tuples on $X$ (notwithstanding that a good analogue of
$\hat{X}^n$ is not known if $\dim X > 1$).

Here we will realize this philosophy as follows.  First, we
construct certain spaces $X<n>$, related to symmetric products, which
we call the {\it very symmetric products} of $X$.  To be precise, $X<n>$
parametrizes the nonempty subsets of $X$ of cardinality $\le n$.  These naturally form
a tower:
\begin{eqnarray*}
X = X <1> \subset X_2 = X <2> \subset X<3> \cdots \subset X<n>  \cdots
\subset X <\infty> = \\
\lim\limits_\to X<n> .
\end{eqnarray*}
Then on $X<\infty>$ we construct a certain complex $J^\cdot = J^\cdot (T_X)$
which we call the {\it Jacobi complex} of $X$:  this is essentially just a
multivariate version of the standard complex used to compute the {\it Lie
algebra homology} of $T_X$. cf. [F] (indeed the latter homology coincides with
the cohomology of $J^\cdot$ along $X<1>$). The subcomplex $F^n J^\cdot =:
J^\cdot_n$ is natually supported on $X<n> \subset X<\infty>$. With this, we
will prove the following

\begin{theorem}
Let $X$ be a compact complex manifold with $H^0(T_X) = 0$ and let $J$ be the
Jacobi complex of $X$. Then

(i) for each n there is a canonical ring structure on
$$R^u_n = {\mathbb C} \oplus H^0(J_n)^* $$
and a canonical flat deformation $X^u_n/R^u_n$,\;these fit together
to form a direct system with limit
$$\hat{X}^u/\hat{R}^u = \lim_n X^u_n/R^u_n\;;$$

(ii) for any artin local ${\mathbb C}-algebra\ R_n$ of exponent n and flat deformation
$X_n/R_n$ of $X$, there is a canonical Kodaira-Spencer ring homomorphism
$$\alpha_n = \alpha_n(X_n/R_n): R^u_n \to R_n$$
and an isomorphism
$$X_n/R_n \stackrel{\sim}{\to} \alpha^*_nX^u_n = X^u_n
\times_{R^u_n} R_n\;;$$

(iii) if $\hat{R}=\lim\limits_\leftarrow R_n$ is a complete local noetherian
${\mathbb C}-algebra$ and $\hat{X} = \lim\limits X_n/R_n$, then $\hat{\alpha} = \lim_n
\alpha_n : \hat{R}^n \to \hat{R}$ exists and $\hat{X}/\hat{R} = \hat{\alpha}
^*( \hat{X}^u/\hat{R}^u)$
\end{theorem}

The result naturally generalizes (cf. Sect. 5). If $g$ is a sheaf of ${\mathbb C}$-
Lie algebras on $X$ and $E$ a $g$-module (both assumed reasonably 'tame'),
let $\bar{g}$ be the unique quotient of $g$ acting faithfully on $E$ and assume
that $H^0(\bar{g})=0$. For
any artin local ${\mathbb C}$-algebra $(R,m)$ of exponent $n$ we may
define a sheaf of groups $G_R$ on $X$ by
$$G_R = \exp (g\otimes m)$$
(i.e. $G_R$ is $g\otimes m$ with multiplication determined by the Campbell-
Hausdorff formula) and similarly $\bar{G}_R = \exp (\bar{g}\otimes m)$.
Then $g$-deformations of $E$ over $R$ are locally trivial deformations with
transitions in $\bar{G}_R$, and are naturally classified by the nonabelian
Cech cohomology set $H^1(X,\bar{G}_R)$. Our construction yields a bijection
$v = v_{R,E}$ between these and a certain subset of ${\mathbb
H}^0(J_n(\bar{g}))\otimes m$ (i.e. the set of 'morphic' elements). For $n\geq 3$
this correspondence is given somewhat indirectly and in particular does not
come from an explicit correspondence on the cocycle level. $v$ apparently
depends on both $E$ and $g$, though its source and target depend only on $\bar{g}$.
It is unknown to the author whether (say for $E$ a faithful $g$-module) $v$ is
independent of $E$, a fortiori whether it can be defined in terms of $g$ alone.
For another, perhaps more 'conceptual' interpretation of our construction of
the universal deformation of $E$, see [R3], Theorem 3.1.

A further generalization, to the
case of a sheaf of differential graded Lie algebras, will be considered in
[R3].

As indicated above, the existence of the universal formal deformation
$\hat{X}^u/\hat{R}^u$ was known before, thanks to the work of
Grothendieck, Schlessinger et al.: our point is its explicit construction and
description.
As for applications and extensions of the method, these
have been, and will be given elsewhere, but a few can be mentioned here.

(i) An analogous deformation theory for deformations of vector
bundles (or more generally locally free sheaves over a fixed
locally C-ringed space) and, as one application, construction of
a symplectic (closed) 2-form on the moduli space, generalising
at the same time constructions of Hitchin (for local systems
on Riemann surfaces) and Mukai (for holomorphic vector bundles
over K3 surfaces)[R5].

(ii) A direct construction of the universal variation of Hodge
structure associated to a compact Kahler manifold and resulting study of
the (local) period map and characterisation of its image (local
Schottky relations), especially for Calabi-Yau manifolds and curves
[R3],[L].

(iii) A theory of semiregularity for submanifolds and embedded, as
well as relative deformations and resulting dimension bounds for
Hilbert schemes and relative deformation spaces [R2][R4].

Higher-order Kodaira-Spencer maps, especially associated to
'geometric' (reduced) families have been independently
defined by Esnault and Viehweg , cf. [EV]; that paper defines higher-
order (additive) Kodaira-Spencer classes, but does not construct
the universal family. Some important antecedents (albeit from a different
viewpoint) are in
the work of Goldman-Millson [GM].

\section{Coalgebra}

The purpose of this section is to characterize the vector space $m^*$ dual
to the maximal ideal of an artin local ${\mathbb C}$-algebra $(R,m)$. While the general
concept of coalgebra is well known, our application in the artin local case
assigns a special role to the m-adic filtration and its dual, the 'order'
filtration,not present in the general case. Consequently, it will be convenient
to give a brief self-contained treatment here.

By an Order-Symbolic (OS) structure of order n we mean  a finite-dimensional
${\mathbb C}$-vector space together with an increasing filtration
$$V^0 = 0\subseteq V^1 \subseteq\cdots \cdots\subseteq V^n=V$$
and mutually compatible 'symbol' or 'comultiplication' maps
$$\sigma^{i,j} : V^i/V^j \to S^2(V^{i-j}),\hskip 30mm j<i .$$
(Sometimes we shall use the same notation to denote the induced map $V^i\to
S^2(V^i)$; actually, a moment's thought shows that $\sigma :=\sigma^{n,1}$
is sufficient to determine the rest) . These are assumed to satisfy
the natural (co)associativity condition that the following diagram should
commute
$$\begin{array}{cccc}
&S^2(V^{n-1})\otimes V^{n-1}&\\
\hskip 50mm\stackrel{\phi}{\nearrow}& &\searrow\hskip 40mm\\
V/V^1 \to S^2V^{n-1} \subset V^{n-1}\otimes V^{n-1}&  & V^{n-1}\otimes V^{n-1}
\otimes V^{n-1}\\
\hskip 50mm\stackrel{\psi}{\searrow}& &\nearrow\hskip 40mm\\
&V^{n-1}\otimes S^2(V^{n-1})&
\end{array}$$
$$\varphi = \sigma^{n-1,1}\otimes id,\hskip 20mm \psi = id \otimes \sigma^{n-1,1}$$

An OS structure V is said to be {\it standard} if $\sigma$ is injective ( hence
 $\sigma^{n,j}$ is injective for all $j < n$). We can now state the basic
result about OS  structures, which relates them with artin local algebras.

\begin{propo}
There is an equivalence of categories between ${\mathcal OS}_n$, the category of OS
structures of order n, and ${\mathcal FR}_n$, the category of commutative artin local
${\mathbb C}$-algebras space $(R,m)$ of exponent n together with a super-m-adic filtration
$(m_i\supseteq (m)^i)$, where standard structures correspond with m-adically
filtered algebras. The correspondence is given by
$$(V,{V^\cdot},\varphi) \mapsto ({\mathbb C}\oplus V^*, V^{i\perp } = (V/V^i)^*, \varphi^*_n)$$
$$(S, m,m_\cdot) \mapsto (m^*, m^\perp_{i-1} = (m/m_{i-1})^*,
{\rm comultiplication}).$$
\end{propo}

\begin{proof}
Basically trivial. Given V etc. define
$$R = {\mathbb C}\oplus V^*, m=V^*, m_i = (V^{i-1})^ \perp = (V/V^{i-1})^*\subset V^*$$
Dualising $\sigma$ yields the multiplication map
$$S^2(m) \to S^2(m/m_n) \stackrel{\sigma^*}{\rightarrow} m_2 \subset m$$

This extends in an obvious way to a commutative associative multiplication
map $S^2R \to R$. By construction, $\sigma^*$ descends to a map
$$S^2(m/m_i) \stackrel{\sigma^{i+1,2 *}}{\longrightarrow} m_2/m_{i+1}$$
hence $m\cdot m_i \subset m_{i+1}$. So inductively $m_i$ is firstly an ideal and
then $m_i \supseteq m^i$ by induction. The rest is similar.
\end{proof}

Thus in particular, to an artin local ${\mathbb C}$-algebra $(R,m)$ of exponent n, we have
a uniquely determined standard OS structure on $T^nR = m^*$, which conversely
determines $(R,m)$. For later use it is convenient to explicate and amplify
the morphism part of the above equivalence.

\begin{corol}
Let $(R,m),(R^\prime, m^\prime)$ be artin local algebras of exponent n. Then
the following are mutually interchangeable:

(i)  a local homomorphism $\eta: R^\prime \to R$;

(ii) an OS morphism $\kappa:T^nR \to T^nR^\prime$;

(iii) a compatible collection of elements
$$v_i \in m^{n+1-i} \otimes T^nR^\prime / T^{n-i}R^\prime$$
such that \hskip 4mm $(id \otimes \sigma)(v_n) =v_n\cdot v_n \in m^2 \otimes S^2(T^nR^\prime)$ \hskip 40mm
(1.1)
\end{corol}

\begin{proof}
Only (iii) may require comment. $v_n$ evidently determines $\kappa$ as well as
$v_1,\cdots \cdots, v_{n-1}$; it is the existence of the latter that ensures that
$\kappa$ is filtration-preserving, while (1.1) makes $\kappa$ compatible with
comultiplication.
\end{proof}

Let us call an element $v \in m_R\otimes T^nR^\prime$ as above {\it morphic}.

\section{Products}
\subsection{Very symmetric products}
Fix a topological space X. For any $n \geq 1$, we denote by $X^n$ and $X_n$
the Cartesian and symmetric products, respectively. The system $(X^n , n \in
{\mathbb N})$ forms essentially a {\it simplicial} configuration ( while the $X_n$'s are
related to one  another in even more complicated ways). On the other hand,
 the system of the n-th order neighborhoods of a point(say on a moduli space),
$ n \in {\mathbb N}$, is simply a {\it  tower}. This indicates that the
'right' spaces of point-configurations to work with in deformation theory
are neither $X^n$ or $X_n$ but a suitable modifications thereof which form
a tower. We now proceed to define these spaces which we call the
{\it very symmetric products} (powers) of X and denote by $X<n>$. A word to the
wise: defining $X<n>$ may appear to be a fastidious bother as (sheaf) cohomology
behaves simply with respect to finite maps; however, it is {\it  complexes}
that we must work with, and even to define the
coboundary maps in appropriate complexes, a certain minimum amount of
'quotienting' must be effected , e.g. it seems that the Jacobi complexes
defined below on $X<n>$ cannot be defined on any natural space strictly
'above' $X<n>$ (and this certainly includes $X_n$).

As a set , we define

$X<n> = X^n / \sim$

$(x_1,\cdots , x_n) \sim (y_1,\cdots ,y_n)$   iff $\{x_1,\cdots,x_n\} = \{y_1,
\cdots, y_n\}$.

Thus $X<n>$ parametrises precisely the nonempty subsets of $X$ of
cardinality $\leq n$. We endow $X<n>$ with the quotient topology
induced from $X^n$. Note that we indeed have a tower of ( closed,
for $X$ separated) embeddings
\begin{eqnarray*}
X = X <1> \subset X_2 = X <2> \subset X<3> \cdots \subset X<n>  \cdots
\subset X <\infty> = \\
\lim\limits_\to X<n> .
\end{eqnarray*}

Alternatively, $X<n>$ may be defined inductively: let
$$X^{n-1} \to X^n \to X_n$$
be a 'diagonal' map, e.g. $(x_1,\cdots,x_{n-1},x_{n-1}) \mapsto \{x_1,\cdots,x_{n-1}\}$
,whose image $D_{n-1}$ is the big 'diagonal' in $X_n$ (and is independent of
the choice of which point gets doubled). Then
\begin{eqnarray*}
X<n> & = & X_n\bigcup\limits_{X^{n-1}} X<n-1> \hskip 20mm (2.1)\\
     & = & X_n\bigcup\limits_{D^{n-1}} X<n-1>
\end{eqnarray*}

(it is easy to see inductively that the natural map $q_n: X_n \to X<n>$ factors
through $D_n$).  Via (2.1), very symmetric products may be defined in more general
settings, e.g. when $X$ is a Grothendieck topology.

It is not hard to see that if $X$ has a structure of (separated) analytic space,
then so does $X<n>$ in a natural way. However, we shall not need this fact.
Rather, the sheaves on $X<n>$ relevant to us will be alternating products of
sheaves induced from $X$, which now proceed to define. Let S be a ring and A a
sheaf of S-modules on $X$. Let $\pi_n: X^n \to X<n>$ be the natural map, and set
$$\tau^n_S (A) = \pi_{n*} ( A\boxtimes_S\cdots\boxtimes_SA)$$
( When S is understood, e.g. $S={\mathbb C}$, we may suppress it). Note that
the symmetric group $\Sigma_n$ acts in a natural way on $\tau^n_S(A)$ and let
$\sigma^n_S(A)$ (resp. $\lambda^n_S(A)$) denote the invariant and antiinvariant
factors. Note that this definition makes sense on the symmetric product $X_n$
already, and may also be extended to mixed ( Schur ) tensors in an obvious way.

When $A$ is replaced by a complex $A^.$ of $S$-modules, these constructions extend
in a natural way to make $\tau^n_S(A^.), \sigma^n_S(A^.), \lambda^n_S(A^.)$ into
complexes; for instance
$$\lambda^2_S(A^.) = \lambda^2_S(A^{even})\oplus \pi_{2*}(A^{even}\boxtimes A^{odd})
\oplus \sigma^2_S(A^{odd}).$$

The cohomology of $\tau_S^n(A)$ can be computed by the K${\rm\ddot{u}}$nneth formula, at least
if $A$ is S-free, i.e.
$$ H^m( \tau_S^n(A)) = [\otimes^n_1 H^\cdot(A)]^m$$
In fact the n-th tensor power of a $\breve{C}$ech complex for A ( with respect
to an acyclic cover of $X$) yields one for $\tau^n_S(A)$. As everything
decompose into  $\pm$ eigenspaces under the action of $\Sigma_n$, analogous
comments apply to $\sigma^n_S(A)$ and $\lambda^n_S(A)$ (one must take into
account the usual sign rules for cup products, e.g. $a\cup b=(-1)^{\deg a\deg b}
b\cup a$). For instance, in the case of principal interest to us, we have $H^0(A)=
0$ and then
$$H^i(X<n>,\lambda^n_S(A))= H^i(X<n>,\tau^n_S(A)) = 0, i<n;$$
$$H^n(X<n>, \lambda^n_s(A)) = S^n_SH^1(A) \ :$$
in fact, the symmetric power of the $\breve{C}$ech complex for A may be used to
compute the cohomology of $\lambda^n_S(A)$.

{\bf Remark} The spaces $X<n>$ have recently appeared in the work of
Beilinson and Drinfeld on 'Chiral Algebras'; I am grateful to V. Ginzburg
for pointing this out.

\subsection{Jacobi complex}
Let ${\mathcal L}^.$ be a sheaf of complex differential graded Lie algebras
(DGLAs) on X. Thus ${\mathcal L}^.$ is a "lie object' in the category of
complexes of ${\mathbb C}$-modules on $X$, which means there is a morphism
$bt:\Lambda^2({\mathcal L}^.)\to {\mathcal L}^.$, whose natural extension as
a derivation of degree -1 on the Grassmann algebra
$\oplus\Lambda^i({\mathcal L}^.)$ satisfies $bt^2=0$; or course 'Grassmann
algebra' and wedge must be understood in the graded sense, compatible with
the gradation on ${\mathcal L}^.$.Note that bt induces a map
$$br:\lambda^2({\mathcal L}^.)\to \Lambda^2({\mathcal L}^.)\to {\mathcal L}^.$$
(i.e. restriction followed by bt).  Now we associate to
${\mathcal L}^.$ a complex $J^\cdot({\mathcal L}^.)$ on $X<\infty>$ called the
Jacobi complex of ${\mathcal L}^.$, as follows. Set
$$J^{-n}({\mathcal L}) = \lambda^n({\mathcal L}^.), \hskip 10mm  n\geq1$$
( where the latter is viewed as a sheaf on $X<\infty>$ via $X<n>\subset X<\infty>$
and $\lambda$ is understood in the graded sense);
the differential $d_n:\lambda^n({\mathcal L}^.) \to \lambda^{(n-1)}({\mathcal L}^.)$
is defined as follows. First, let $alt:\tau^2({\mathcal L}^.)\to \lambda^2({\mathcal L}
^.)$ be the alternation or skew-symmetrization map, where $\lambda^2({\mathcal L}^.)$
is viewed as a complex on  $X<2>$ via the diagonal embedding $X\to X<2>$, and set
$$a = \pi_{n*}(id \boxtimes alt):\pi_{n*}(\boxtimes^n {\mathcal L}^.) \to
\pi_{n*}(\boxtimes^{n-2}{\mathcal L}^.\boxtimes\lambda^2({\mathcal L}^.)).$$
Next, note that $\pi_{n*}(\boxtimes^{n-2}{\mathcal L}^.\boxtimes\lambda^2({\mathcal L}
^.))$ as defined above coincides with $\pi_{(n-1)*}(\boxtimes^{n-2}{\mathcal L}^.
\boxtimes \lambda^2({\mathcal L}^.))$ and set
$$b = \pi_{(n-1)*}(id\boxtimes br):\pi_{(n-1)*}(\boxtimes^{n-2}{\mathcal L}^.
\boxtimes\lambda^2({\mathcal L}^.))\to \pi_{(n-1)*}(\boxtimes^{n-2}{\mathcal L}^.
\boxtimes {\mathcal L}^.)=\tau^{n-1}({\mathcal L}^.).$$
Finally let $p:\tau^{n-1}({\mathcal L}^.) \to \lambda^{n-1}({\mathcal
L}^.)$ be the natural alternation map and $i:\lambda^n({\mathcal L}^.)\to
\tau^n({\mathcal L}^.)$ the inclusion. Then define
$$d_n = p\circ b\circ a\circ i .$$

More explicitly,
$$d_n(t_1\times\cdots\times t_n) = \frac{1}{n!}\sum\limits_{\sigma\in \Sigma_n}
\sgn(\sigma)[t_{\sigma(1)}, t_{\sigma(2)}]\times t_{\sigma(3)}\times
\cdots\times t_{\sigma(n)}$$

The Jacobi identity for ${\mathcal L^.}$ ensures that $J({\mathcal L^.})$ is a
complex. Put $J^\cdot_n({\mathcal L}^.) = J^{\geq -n}({\mathcal L}^.)$, which may
be viewed as complex
on $X<n>$.

Now note that,\; viewing \;$J_n({\mathcal L}^.)/J_1({\mathcal L}^.)\;
=\; J^{-n\leq \cdot \leq -2}({\mathcal L}^.)$
as a complex on \\$X<2n-2>$,  it forms a subcomplex of $\pi_{n-1,n-1*}(Sym^2(J
_{n-1}({\mathcal L})^.))$, where $\pi_{n-1,n-1} : X<n-1><2> \to X<2n-2>$ is the
natural map. This gives rise, e.g. to a map
$$ \sigma^n: {\Bbb H}^0(J_n({\mathcal L}^.))/ {\Bbb H}^0(J_1({\mathcal L}^.))
\to {\Bbb H}^0((J_n/J_1)({\mathcal L}^.)
) \to S^2{\Bbb H}^0(J_{n-1}({\mathcal L}^.))$$
which we call the symbol map associated to ${\mathcal L}^.$, it is not hard to
see that with this $V^n({\mathcal L}^.) = {\Bbb H}^0(J_n({\mathcal L}^.))$ forms an
OS structure, which
is standard provided ${\Bbb H}^{\leq 0}({\mathcal L}^.) = 0$. By Section 1 then we obtain
an inverse system
of artin local algebras
$$R_n({\mathcal L}^.) = {\mathbb C}\oplus V^n({\mathcal L}^.)^*$$
and their limit $\hat{R}({\mathcal L}^.)$ which might be called the deformation
ring associated to ${\mathcal L}^.$.

In particular, if $X$ is a compact complex manifold, its tangent sheaf $T = T_X$
forms a  Lie algebra under Lie bracket of vector fields, and we denote the
associated Jacobi complexes by $J_{n,X}$ or $J_n$, the corresponding OS structure
by $V^n_X$ or $V^n$, and the corresponding ring by $R^u_{n,X}$ or $R^u_n$,
As we shall see, when $H^0(T) = 0$ the latter turns out to be the base ring of
the n-universal deformation of $X$.
\subsection{Obstructions}
Assume $H^0({\mathcal L})=0$. Note that the long cohomology
sequence associated to $$0\to J_{n-1}({\mathcal L})\to
J_n({\mathcal L})\to \lambda^n({\mathcal L})[n]\to 0, n\geq 2,$$
gives rise to a 'big obstruction' map $$Ob_n:Sym^nH^1({\mathcal
L})\to {\mathbb H}^1(J_{n-1}({\mathcal L})).$$ Let $K^n=\ker
(Ob_n)$, so that we have an exact sequence $$0\to V^{n-1}\to
V^n\to K^n\to 0.$$ Then, using the 'comultiplicative' structure on
$J_n({\mathcal L})$ as above it is easy to see that $Ob_n$ factors
through a map, denoted $ob_n$, called the 'small obstruction map
$$ob_n:K_{n-1}. H^1({\mathcal L})\to H^2({\mathcal L}),$$ where
$K^{n-1}.H^1({\mathcal L})$ denotes the intersection of
$Sym^nH^1({\mathcal L})$ and $K^{n-1}\otimes H^1({\mathcal L})$
considered as subspaces of $\otimes^nH^1({\mathcal L})$, and we
have $K^n=\ker ob_n$ as well. Thus $K^n$ may be described
inductively, starting with $K^1=H^1({\mathcal L}), ob_1=0$.

Similar comments apply if ${\mathcal L}$ is replaced by a dgla
${\mathcal L}^.$ with ${\mathbb H}^{\leq 0}({\mathcal L}^.)=0$.

\section{Second order}
For $n = 1$, Theorem 0.1 reduces to standard first -order Kodaira-Spencer
deformation theory. Before taking up the general n-th order case in the next section,
we consider here the second-order case , which is relatively simple but already
illustrates some of the ideas. Thus let us fix an artin local ${\Bbb C}$-algebra$(
R_2,m_2)$ with reduction $(R_1,m_1) = (R_2/m^2_2, m_2/m^2_2)$ as well as an acyclic
(say polydisc) open cover $(U_\alpha)$ of $X$, to be used in computing $\breve{C}$ech
cohomology. To a flat deformation
$$ X_2/R_2 = Spec({\mathcal O}_2)$$
we seek to associate a Kodaira-Spencer homomorphism
$$\alpha_2 = \alpha_2(X_2/R_2): R^u_2 \to R_2$$
or equivalently (cf. Section 2)  a morphic element
$$v_2=v_2( X_2/R_2) \in m_2 \otimes H^0(J_2)$$
which is to be described by a hypercocycle

$v_2 = (u, \frac{1}{2}u^2) \in \breve{C}^1(T)\otimes m_2 \oplus S^2\breve{C}
^1(T)\otimes m^2_2 \subset \breve{C}^0(J_2)\otimes m_2 \hskip 18mm$ (3.1)\\
where $u=(u_{\alpha\beta})$ is required to be a lifting of
$$ v_1 = (v_{1\alpha\beta})\in \breve{Z}^1(T)\otimes m_1,$$
a cocycle representing ${\mathcal O}_1 ={\mathcal O}_2\otimes_{R_2} R_1$
(where ${\mathcal O}_2$ is the structure sheaf of $X_2$), and $u^2$ means exterior
cup product in the cochain sense,i.e
$$(u^2)_{\alpha\beta\gamma} = u_{\alpha\beta}\times u_{\beta\gamma} \in
S^2\breve{C}^1(T)\otimes m^2_2 \subset \breve{C}^2(\lambda^2T)\otimes m^2_2$$
Note that the particular form of (3.1) makes the morphicity of $v_2$ automatic
provided it is a hypercocycle, which means explicitely

$$ -\frac{1}{2}[u_{\alpha\beta}, u_{\beta\gamma}] = u_{\alpha\beta}+u_{\beta\gamma}+
u_{\gamma\alpha} = \delta(u) \hskip 30mm \text{(3.2)}$$
$$ \delta = \breve{C}\text{ech\ coboundary}$$
Note that the LHS of (3.2) depends only on the reduction $v_{1\alpha\beta}$
of $u_{\alpha\beta}$ mod $m^2_2$

Now to define $(u)$ we proceed as follows. As ${\mathcal O}_2/R_2 $ in a flat deformation
of ${\mathcal O}$, it is locally trivial hence we have isomorphisms of $R_2-$algebras
$$\psi_\alpha : {\mathcal O}_2(U_\alpha) \to {\mathcal O}(U_\alpha)\otimes_{\mathbb C}R_2$$
which give rise to a gluing cocyle given by
$$ D^2_{\alpha\beta} =\psi_\alpha\circ\psi^{-1}_\beta \in Aut_{R_2}({\mathcal O}(U_\alpha
\cap U_\beta)\otimes R_2)$$
which reduces mod $m^2_2$ to
$$D^1_{\alpha\beta} = I + v_{1\alpha\beta} \in Aut_{R_1}({\mathcal O}(U_\alpha\cap U_\beta)
\otimes R_1),$$
 a gluing cocycle defining ${\mathcal O}_1$.

Now it is easy to see that $D^2_{\alpha\beta}$ is uniquely expressible in the
form
\begin{eqnarray*}
D^2_{\alpha\beta}=&\exp(u_{\alpha\beta}) & \hskip 30mm (3.3)\\
   =& I + u_{\alpha\beta} +\frac{1}{2}u^2_{\alpha\beta}, &\hskip 6mm u_{\alpha\beta} \in
m_2\otimes T(U_\alpha\cap U_\beta):
\end{eqnarray*}
indeed starting with an arbitrary lift $u^\prime_{\alpha\beta}$ of $v_{1\alpha\beta}$
to $m_2\otimes T(U_a\cap U_\beta)$, $\exp(u^\prime_{\alpha\beta})$ and $D^2_{\alpha\beta}$
are $R_1$-algebra homomorphisms which agree mod $m^2_2$, hence differ by an $m^2_2$-
valued derivation $t_{\alpha\beta}$ and we may set $u_{\alpha\beta}= u^\prime_{\alpha\beta}
+ t_{\alpha\beta}$. Now we simply plug (3.3) into the cocycle equation for $D^2$:
$$D^2_{\alpha\beta}D^2_{\beta\gamma} = D^2_{\alpha\gamma} \hskip 50mm (3.4)$$
which becomes,
\begin{eqnarray*}
I +u_{\alpha\beta} +u_{\beta\gamma} +\frac{1}{2}u^2_{\alpha\beta} +u_{\alpha\beta}
u_{\beta\gamma}+\frac{1}{2}u^2_{\beta\gamma} &=& I +u_{\alpha\gamma} +\frac{1}{2}(u_
{\alpha\gamma})^2\\
&=&I + u_{\alpha\gamma} +\frac{1}{2}(u^2_{\alpha\beta}+u_{\alpha\beta}u_{\beta\gamma}
+u_{\beta\gamma}u_{\alpha\beta}+u^2_{\beta\gamma})
\end{eqnarray*}
as $(u_{\alpha\beta})$ is a cocycle mod $m^2_2$. This is obviously equivalent
to (3.2). Thus $v_2$ is a hypercocycle , as claimed.

Now the foregoing argument can essentially be read backwards given a morphic
element
$$v_2 \in m_2\otimes H^0(J_2), $$
choose a representative for $v_2$ of the form
$$((u_{\alpha\beta}), (u^\prime_{\alpha\beta\gamma})) \in \breve{C}^\prime(T)\otimes m_2
\oplus S^2\breve{C}^\prime(T)\otimes m^2_2 \subset \breve{Z}^0(J_2),$$
where $(u_{\alpha\beta})$ is a lifting of $(v_{1\alpha\beta})$; thus compatibility
with comultiplication yields that $v_2$ may also be represented by
$$((u_{\alpha\beta}),\frac{1}{2}(u_{\alpha\beta})^2).$$
Then simply setting $D^2_{\alpha\beta}= \exp(u_{\alpha\beta})$, the cocycle
condition (3.4) follows from the hypercocycle condition (3.2), so that $(D^2_{\alpha\beta})$
yields a locally trivial flat deformation $X_2/R_2=Spec {\mathcal O}_2$, which we denote by
$\Phi_2(\alpha_2)$ ( though it is yet to be established that this is independent of
choices).

This construction applies in particular to the identity map $R^u_2 \to R^u_2$,
 thus yielding a flat deformation over $R^u_2$ which we call an
{\it  universal second order deformation} and denote by $X^u_2=Spec(
{\mathcal O}^u_2)$. It is moreover clear by construction that $\Phi_2(\alpha_2) =
\alpha^*_2(X^\kappa_2/R^n_2)$ for any $\alpha_2: R^u_2\to R_2$ and also
that for any second-order deformation $X_2/R_2$,
$$X_2/R_2 \approx \alpha_2(X_2/R_2)^*(X^u_2/R^u_2) \approx \Phi_2(
\alpha_2(X_2/R_2)).$$
Similarly,
$$\alpha _2(\Phi_2(\beta))= \beta.$$
Thus $\alpha_2$ and $\Phi_2$ establish mutually inverse correspondences, albeit
on the cocycle level. What has to be established is that this correspondence
descends to cohomology , i.e. non-abelian cohomology of Aut-cocycles and
hypercohomology respectively. In one direction, consider two cohomologous Aut-cocycles
$$D^2_{\alpha\beta} \sim  D^{2\prime}_{\alpha\beta} = A_\beta D^2_{\alpha\beta}
A^{-1}_\alpha$$
$A_\alpha \in Aut_{R_2}({\mathcal O}(U_\alpha)\otimes R_2)$, as above uniquely expressible
in the form $\exp(w_\alpha), w_\alpha \in m_2\otimes T(U_\alpha)$. Thus
\begin{eqnarray*}
D^{2\prime}_{\alpha\beta} &=& (I +w_\beta +\frac{1}{2}w^2_\beta)(I + u_{\alpha\beta}
+\frac{1}{2}u^2_{\alpha\beta})(I - w_\alpha + \frac{1}{2}w^2_\alpha)\\
&=&I + (u_{\alpha\beta} +w_\beta -w_\alpha +\frac{1}{2}[w_\beta- w_\alpha, u_{
\alpha\beta}] +\frac{1}{2}[w_\alpha, w_\beta]) +\frac{1}{2}(u_{\alpha\beta}+w_
\beta -w_\alpha)^2\\
&=& \exp(u_{\alpha\beta}+ w_{\beta} -w_{\alpha} + \frac{1}{2}[w_\beta -w_\alpha, u_{\alpha\beta}]
+ \frac{1}{2}[w_\alpha, w_\beta])\\
&=&: \exp(u^\prime_{\alpha\beta})
\end{eqnarray*}
Then $v^\prime_2 = v_2(D^{2\prime}) = (u^\prime, \frac{1}{2}(u^\prime)^2)$ is
cohomologous to $v_2$ because
$$v^\prime_2-v_2 = \partial((w_\alpha),\frac{1}{2}(w_\alpha\times u_{\alpha
\beta})+ \frac{1}{2}(w_\alpha\times w_\beta))$$
where $\partial = \delta \pm b$ is the differential of the $\breve{C}$ech bicomplex
of $\breve{C}(J_2)$. Conversely, supposing $v_2 = (u, \frac{1}{2}u^2), v^\prime
_2 = (u^\prime, \frac{1}{2}u^{\prime2})$ are cohomologous,
$$v^\prime_2 -v_2 = \partial((w_\alpha),(t_{\alpha\beta}))\;.$$
Now as $ H^0(T) = 0, \delta(t) = \frac{1}{2}(u^\prime)^2 -\frac{1}{2}u^2$
determines $(t)$ up to adding a $\breve{C}$ech coboundary $s_\alpha- s_\beta$ and,
using $b\delta = \pm \delta b$ this may be absorbed into $(w_\alpha)$. Thus we may assume
$$t_{\alpha\beta} = \frac{1}{2}w_\alpha\times u_{\alpha\beta} +\frac{1}{2}
w_\alpha\times w_\beta,$$
so that $(D^2_{\alpha\beta}= \exp(u_{\alpha\beta}))$ and $(D^{2\prime}_{\alpha\beta}=
\exp(u^\prime_{\alpha\beta}))$ are cohomologous as above. This finally completes
the proof of Theorem 0.1 for n=2.

\section{$n$-th order}
We now complete the proof of Theorem 0.1 in the general n-th order case,
following in part
the pattern of the case n=2 and using induction.
However the argument becomes a bit
more involved and less direct. Fix an artin local ${\mathbb
C}$-algebra $(R_n, m_n)$ of exponent n, with reduction $(R_{n-1}, m_{n-1})$, etc,
and an acyclic open cover $(U_\alpha)$. The main point is to associate a morphic
hypercocycle
$$v_n = v_n({\mathcal O}_n/R_n) \in m_n\otimes \breve{Z}^0(J_n)\;,$$
hence a Kodaira-Spencer homomorphism $\alpha_n(O_n/R_n)$ etc- to an $R_n$-flat
deformation
${\mathcal O}_n={\mathcal O}_{X_n}$ of ${\mathcal O}$.
As before we seek $v_n$ of the form
$$v_n = \epsilon(u_n) :=(u_n, \frac{1}{2}(u_n)^2,\cdots,\frac{1}{n!}(u_n)^n)$$
for some cochain $u_n = (u_{n\alpha\beta}) \in \breve{C}^1(T)\otimes m_n$
which is a lift of $u_{n-1} \in \breve{C}^1(T)\otimes m_n$ analogously
defining $v_{n-1}$. To this end we start with isomorphisms of algebras
$$\psi^n_\alpha: {\mathcal O}_n(U_\alpha) \stackrel{\sim}{\to} {\mathcal O}(U_\alpha)\otimes R_n$$
which yield a gluing cocycle by
$$D^n_{\alpha\beta} = \psi^n_{\beta}(\psi^n_\alpha)^{-1} \in Aut_{R_n}({\mathcal O}(U_\alpha
\cap U_\beta)\otimes R_n), \hskip 20mm (4.1)$$
which as above we express in the form
$$D^n_{\alpha\beta} = \exp(t_{n\alpha\beta}), \hskip 58mm (4.2)$$
This can be done because, assuming inductively  that (4.2) holds for $n-1$ and
letting $t^\prime_n$ be an arbitrary lift of $t_{n-1}$ and $t_n = t^\prime_n+
\eta_n, \eta_n \in \breve{C}^1(T)\otimes m^n_n,$ (4.2) can be rewritten as
$$D^n_{\alpha\beta} = \exp(t^\prime_{n\alpha\beta}) + \eta_n.$$
which can clearly be uniquely solved for $\eta_n$.

Now before proceeding with the definition of $u=u_n$ we will
consider a Dolbeault analogue, both for its own interest and as
motivation for the Cech construction to follow.
Consider the DGLA sheaf
$$ g^. = ({\mathcal A}^{0,.}(T),\bar\partial , [, ] ) $$
($\Gamma (g^.)$ is sometimes called the Frohlicher-Nijenhuis
algebra); as $g^.$ is a soft resolution of $T$, $J_n(g^.)$ is
a soft resolution of $J_n(T)$ which may be used to compute
${\Bbb H}^0(J_n(T))$. As $g^0$ is soft it is easy to see
that, up to shrinking our cover $(U_{\alpha})$ we may assume
$$D^n_{\alpha\beta} = exp(s_{\alpha})exp(-s_{\beta})$$
$$s_{\alpha}\in g^0(U_{\alpha})\otimes m_n. $$

Put another way, we may view $\psi ^n_{\alpha}$ above as a
holomorphic local
trivialisation
$$
U^n_{\alpha} \simeq U_{\alpha}\times {\rm Spec}(R_n)$$
$U^n_{\alpha} = $ open subset of $X_n$ corresponding to
$U_{\alpha}$; on the other hand there is a global '$C^{\infty}$
trivialisation' $C: X_n \to X\times {\rm Spec}(R_n)$,
and we may set
$$
\exp (s_{\alpha}) = C\circ (\psi^n_{\alpha})^{-1} \hskip 20mm (4.3)$$

Now note that $\bar\partial$ extends formally as a derivation on
the universal enveloping algebra $U(g^.)$ and we set
$$ \phi _{\alpha} = exp(-s_{\alpha})\bar\partial exp(s_{\alpha})=
D(ad(s_{\alpha}))(\bar\partial s_{\alpha}) \hskip 20mm (4.4)$$
where $D$ is the function
$$ D(x) = \frac{exp(x)-1}{x} = \sum_{i=0}^{\infty} {\frac{x^i}{(i+1)!}}. $$
Note that
$$
0=\bar\partial D^n_{\alpha\beta}=\bar\partial exp(s_{\alpha}) exp(-s_{\beta})
+exp(s_{\alpha})\bar\partial exp(exp(-s_{\beta}),$$
hence
$$
exp(-s_{\alpha})\bar\partial exp(s_{\alpha})= -\bar\partial
exp(-s_{\beta})exp(s_{\beta}) ;$$
since moreover $\bar\partial (exp(-s_{\beta})exp(s_{\beta}) = 0$
we have similarly
$$
-\bar\partial exp(-s_{\beta})exp(s_{\beta})=exp(-s_{\beta})
\bar\partial exp(s_{\beta}),\hskip 20mm (4.5)$$
which means precisely that the $\phi_{\alpha}$ glue together
to a global section
$$
\phi \in \Gamma (g^1) = A^{0,1}(T)\otimes m_n. $$
Next, note using (4.4) that
$$
\bar\partial \phi _{\alpha} = \bar\partial exp(-s_{\alpha})
\bar\partial exp(s_{\alpha}) = \bar\partial exp(-s_{\alpha})
exp(s_{\alpha})exp(-s_{\alpha}) \bar\partial exp(s_{\alpha})
= - \phi_{\alpha} \phi_{\alpha} ;$$
recalling that for odd-degree elements $\phi ,\psi \in g^. ,
[\phi ,\psi ] = \phi .\psi + \psi .\phi,$ we conclude that
the integrability equation
$$
\bar\partial\phi =\frac{-1}{2} [\phi ,\phi] \hskip 20mm (4.6)$$
is satisfied, and consequently
$$
\epsilon(\phi)=(\phi , \frac{1}{2}\phi\times\phi ,...,
\frac{1}{n!}\phi\times ...\times\phi )\in \Gamma(J_n(g^.)))
\otimes m_n $$
is a hypercocycle, which may be used to define a Dolbeault
analogue of $v_n$ (automatically morphic, due to the 'exponential'
nature of $\epsilon$).

By way of interpretation, note that, as operators,
$$
\bar\partial (exp(s_{\alpha})) = [\bar\partial , exp(s_{\alpha})],$$
therefore clearly
$$
\phi . = exp(-s.)\bar\partial exp(s.) - \bar\partial .\hskip 20mm (4.6)$$
What (4.6) means is this: recall the map $C$ above which yields
a $C^{\infty}$
trivialisation of the deformation $X_n/R_n$ and in particular
bundle isomorphisms
$$
{\mathcal A}^{0,.}(X)\otimes R_n \simeq {\mathcal A}^{0,.}(X_n/R_n)
$$
under which the canonical Dolbeault operator $\bar\partial_n$
on the RHS corresponds on the LHS precisely to
$\bar\partial_0\otimes 1 +\phi.$ The integrability equation (4.5)
reads, on the operator level
$$
\bar\partial\phi + \phi\bar\partial = \phi\phi, \bar\partial:=
\bar\partial_0\otimes 1, $$
i.e. is equivalent to $\bar\partial _n^2 = 0$.

Given this, it is now clear how to go backwards. Given
$\phi\in A^{0,1}(T)\otimes m_n$ we may define an operator
$d_n$ on ${\tilde\mathcal A}^{0,.}_n:={\mathcal A}^{0,.}(X)\otimes R_n$
by
$$
d_n = \bar\partial + \phi,$$
and the integrability equation (4.5) guarantees that $({\tilde\mathcal A}
^{0,.}_n,d_n)$ is a complex; by semicontinuity, this complex is
clearly exact in positive degrees (because ${\tilde\mathcal A}^{0,.}_n
\otimes \Bbb C$ is) and we may define
$$ {\mathcal O}_n = \ker (d_n, {\tilde\mathcal A}^{0,0}_n) . $$
As $d_n$ is an $R_n-$ linear derivation, ${\mathcal O}_n$ is a sheaf
of $R_n$- algebras. That ${\mathcal O}_n$ is $R_n$-flat is a
consequence of the following easy observation.
\begin{lemma} Let $R$ be an artin local ring with residue field $k$ ,
$M$ an $R$-module and $M\to N^.$ a flat resolution such that
$M\otimes k\to N^.\otimes k$ is also a resolution.
Then $M$ is flat.
\end{lemma}
\begin{proof} Our assumption implies that $Tor_i(M,k)=0,i>0$.
Now if $P$ is any finite $R$-module then $P$ admits a composition
series with factors isomorphic to $k$, hence $Tor_i(M,P)=0,i>0$.
Finally any $R$-module $Q$ is a direct limit of its finite submodules
and $Tor$ commutes with direct limits, hence $Tor_i(M,Q)=0$,
so $M$ is flat.
\end{proof}

While the above is sufficient for a Dolbeault proof of Theorem 0.1,
it seems desirable to have a translation into the ${\breve C}$ech
language. To this end, we replace $g^.$ by the ${\breve C}$ech
complex $C^.(T)$ which, together with the ${\breve C}$ech differential
$\delta$ and the natural bracket [,] forms a DGLA. By analogy with
$\phi$, we set
$$ u.= u_n = \delta (s.) = exp(-s.)\delta (exp(s.))
= -\delta (exp(-s.))exp(s.). $$

As $C$ is globally defined (albeit nonholomorphic), it commutes
with $\delta$ hence by (4.3)
$$ u_{n\alpha}= (\psi^n_{\alpha})^{-1}\delta\psi^n_{\alpha}$$
so actually $u_n\in {\breve C}^1(T)$, i.e. $\bar\partial (u_n)=0$.
In particular,
$$
\bar\partial (exp(-s.)\delta (exp(s.)) =
-exp(-s.)\bar\partial \delta (exp(s.) . $$
On the other hand  $\delta (\phi .) =0$ yields
$$\delta (exp(-s.))\bar\partial(exp(s.) = -exp((-s.)\delta \bar\partial
(exp(s.) .$$
As $\delta$ and $\bar\partial$ commute it follows that
$$
\delta (exp(-s.)\bar\partial ((exp(s.)
= \bar\partial ((exp(-s.))\delta (exp(s.).$$
Hence
$$
u.\phi . = exp(-s.)\delta (exp(s.)exp(-s.)\bar\partial (exp(s.)
=-\delta (exp(-s.)\bar\partial (exp(s.))$$
$$
=-\bar\partial (exp(-s.))\delta (exp(s.)) = \phi .u. ,$$
i.e. $u. , \phi .$ commute:
$$
[u. , \phi .] = 0 \hskip 25mm (4.8) $$
Now as above we have formally that
$$ \delta (u.) = \delta (exp (-s.))\delta (exp(s.) = \frac {-1}{2}[u.,u.],$$
and therefore $v_n=\epsilon (u_n) \in {\breve C}^0 (J_n(T))$ is a morphic
hypercocycle,which may be used to define the required Kodaira-Spencer
homomorphism $\alpha _n(X_n/R_n): R_n^u \to R_n.$

The interpretation of $u_n$ is analogous to that of $\phi$: i.e.
the operator
$$
\delta + u_n : {\breve C}^.({\mathcal O})\otimes R_n \to {\breve C}^{.+1}
({\mathcal O})\otimes R_n $$
corresponds to the coboundary operator on ${\breve C}^.({\mathcal O}_n)$
under the local trivialisation $(\psi^n_{\alpha})$ above.
Thus to reverse this construction we may proceed analogously as in the
Dolbeault case. Firstly we represent a morphic element
$v_n \in {\Bbb H}^0(J_n)\otimes m_n$ in the form
$$
v_n =\epsilon (u_n),  u_n \in  {\breve C}^1(T)\otimes m_n \hskip25mm (4.9)$$
where $u.=u_n$ satisfies the ${\breve C}$ech integrability equation
$$\delta (u.) = \frac{-1}{2}[u.,u.]. \hskip 30mm (4.10)$$
Now thanks to (4.9), the deformed coboundary operator
$$
\delta ' = \delta + u_n : {\breve C}^.({\mathcal O})\otimes R_n \to
{\breve C}^{.+1}({\mathcal O})\otimes R_n $$
satisfies $(\delta ')^2 =0$, thus making
$({\breve C}^.({\mathcal O})\otimes R_n,
\delta ')$ as well as its sheafy version
$({\breve\mathcal C}^.({\mathcal O})
\otimes R_n,\delta ')$ into complexes where the latter is exact in positive
degrees. Hence as before,
$$
{\mathcal O}_n = \ker (\delta ',{\breve\mathcal C}^0({\mathcal O})\otimes R_n)$$
is a sheaf of flat $R_n$-algebras yielding a flat deformation
$$
\Phi _n(v_n) = X_n/R_n = Specan ({\mathcal O}_n).$$

It is worth noting that the ${\breve C}$ech construction yields the same
deformation as the Dolbeault one: this follows easily from (4.8).
Also, it is clear from the construction that either $\phi$ or $u_n$
determines the deformation $X_n/R_n$ up to isomorphism.

Now we may easily complete the proof as in Sect.3. First, taking the element
$v_n$ corresponding to the identity on ${\Bbb H}^0(J_n)$, we obtain a
corresponding deformation $X_n^u/R_n^u$. Next, given any $X_n/R_n$,
with corresponding $u_n, \alpha _n$, it is clear that
$$
\alpha _n = \alpha _n(\alpha _n^*(X_n^u/R_n^u).$$
Since a deformation is determined by its $\alpha _n$ it follows that
$$
X_n/R_n \simeq \alpha _n^*(X_n^u/R_n^u).$$
Thus $X_n^u/R_n^u$ is $n$-universal. Finally it is clear by construction
that for different $n$ these are mutually compatible so the limit
${\hat X}_n^u/{\hat R}_n^u$ exists and is formally universal, completing the
proof of Theorem 0.1.

{\bf Remark} See [R3] for an 'interpretation' of the construction of
$X^u_n$.

\section{Generalizations}
Let $g$ be a sheaf of ${\mathbb C}$- Lie algebras on $X$
with $H^0(g)=0$ , $E$
a $g$-module. Replacing $g$ by its unique quotient acting
faithfully on $E$, we may assume $E$ is faithful.
We also assume $X, g, E$ are reasonably tame so
cohomology can be computed by ${\breve C}$ech
complexes.
We further assume $g$ and $E$ admit compatible soft resolutions
$g^.,E^.$ where $g^.$ is a DGLA acting on $E^.$.
Typically, $E$ will have some additional
structure and $g$ will coincide with the full Lie algebra of infinitesimal
automorphisms of the structure: e.g. when $E$ is a ring g may be the algebra of
internal derivations. For any artin local ${\mathbb C}$-algebra $(R,m))$,
we have a
Lie group sheaf $Aut^\circ_R(E\otimes R)$ of
$R$-linear automorphisms of $E\otimes R$
which act as the identity on $E=(E\otimes R)\otimes_R {\mathbb C}$, and we assume
given a Lie subgroup sheaf
$$G_R \subset Aut^\circ_R(E\otimes R)$$
with Lie algebra $g\otimes m$, which coincides - by definition if you will -with
the subgroup of structure-preserving automorphisms in $Aut^\circ_R(E\otimes R)$.
Then the above constructions, being essentially formal in nature, carry over to
this setting essentially verbatim, yielding n-universal deformations $E^u_n
/R^u_n, \ n\geq 1$, and a formally universal deformation $\hat{E}^u/\hat{R}^u$.

{\it  Examples} (cf. [R5])

$5.a$. E is a simple locally free finite-rank ${\mathcal O}_X$-module and g is the algebra of
all traceless ${\mathcal O}_X$-linear endomorphisms of E. The deformation obtained is the
usual universal deformation of E as ${\mathcal O}_X$-module.

{\it  Subexamples}

$5.a_1$. ${\mathcal O}_X$ is the ring of locally constant functions on the
topological space X assumed 'nice', e.g. a manifold.
In this case E is a local system ( i.e. a $\pi_1$ representation), and we obtain
its universal deformation as such.

$5.a_2$. ${\mathcal O}_X$ is the sheaf of holomorphic functions on a complex manifold ( or
regular functions on a proper ${\mathbb C}$-scheme). In this case E is a (holomorphic)
vector bundle and we obtain its universal formal deformation as such.

$5.b.$ Let $Y\subset X$ be an embedding of compact complex manifolds, $g= T_{X/Y}$
the algebra of vector fields on $X$ tangent to $Y$ along $Y$, which may be identified
with the algebra of infinitesimal automorphisms( i.e. internal derivations) of
${\mathcal O}_X$ preserving the subsheaf ${\mathcal I}_y$. Assuming $H^0(T_{X/Y})= 0$, we obtain the
universal deformation of the pair $(X,Y)$.

The case of general holomorphic map $f: Y \to X$ may be treated in a similar way
using the algebra $T_f$ (cf. [R2]); in fact it is sufficient for many purposes to
replace f by the embedding of its graph in $Y\times X$(cf.[R4]). On the other hand the
case of deformations of $Y\to X$ with $X$ fixed requires the DGLA formalism and
will be taken up in [R3].

\end{document}